\documentclass[12pt]{amsart}
\usepackage[active]{srcltx}
\usepackage{a4wide}
\usepackage{amssymb}
\usepackage[T1]{fontenc}
\usepackage[utf8]{inputenc}
\usepackage{fourier}
\usepackage{hyperref}
\usepackage{array}
\usepackage{makecell}
\usepackage{float}

\def\D{{\mathbb D}}
\def\C{{\mathbb C}}
\def\Z{{\mathbb Z}}
\def\cb{{\mathcal B}}
\def\cd{{\mathcal D}} 

\newcommand\norm[1]{\left\lVert#1\right\rVert}
\newtheorem{prop}{\sc Proposition}
\newtheorem{lem}[prop]{\sc Lemma}
\newtheorem{thm}[prop]{\sc Theorem}
\newtheorem{cor}[prop]{\sc Corollary}
\newtheorem{rk}{\sc Remark}

    \title{Truncated Dirichlet estimates for Bloch functions}

    \author[I.~Efraimidis]{Iason Efraimidis}
    \address{Departamento de Matemáticas, Universidad Autónoma de Madrid, 28049 Madrid, Spain}
    \email{iason.efraimidis@uam.es}

    \author[A.~Llinares]{Adrián Llinares}
    \address{Departamento de Matemáticas, Universidad Autónoma de Madrid, 28049 Madrid, Spain}
    \email{adrian.llinares@uam.es}

    \author[A.~Mas]{Alejandro Mas}
    \address{Departamento de Matemáticas, Universidad de Alicante, 03690 Sant Vicent del Raspeig, Spain}
    \email{a.mas@ua.es}

    \date{\today}



\begin{document}

    \begin{abstract}
        We study the truncated Dirichlet functional
            \[
                \Phi_n (f) = \sum_{k = 1}^n k |a_k|^2, \quad f(z) = \sum_{k = 0}^\infty a_k z^k,
            \]
            when $f$ belongs to the Bloch space. If $n \geq 27$, we find a normalized Bloch function $f_n$ such that
            \[
                \Phi_n (f_n) > n B_n^2,
            \]
            where $B_n = \max \{ |a_n| : \| f \|_{\mathcal{B}} = 1 \}$, thus proving that an inequality by Kayumov and Wirths cannot be extended to large values of $n$. Nevertheless, we show that the inequality $\Phi_n(f) \leq n B_n^2$ holds for all normalized functions when $n \leq 9$.
    \end{abstract}
    
    \maketitle

    \section{Introduction}

        Let $\D$ be the open unit disc of the complex plane $\C$. The \textit{Bloch space} $\cb$ is defined as the space of all analytic functions in $\D$ such that
            \begin{equation} \label{eqn:BlochSeminorm}
                \rho_\cb (f) = \sup_{|z| < 1} \{ |f'(z)| (1 - |z|^2) \} < \infty.
            \end{equation}
            Note that $\rho_\cb$ is a seminorm, and it is well-known that $\cb$ becomes a Banach space when endowed with the norm $\| f \|_\cb = |f(0)| + \rho_\cb (f)$. An important property of $\cb$ is that it is Möbius invariant, \emph{i.e.}, $\rho_\cb (f) = \rho_\cb (f \circ \varphi)$ for every disk automorphism $\varphi$. Actually, as shown by Rubel and Timoney \cite{RT}, it is the largest Banach space of holomorphic functions (admitting a non-trivial continuous linear functional) with this property.
        
        Another Möbius invariant space is the \emph{Dirichlet space} $\cd$. We say that an analytic function $f$ belongs to $\cd$ if
            \[
                \rho_\cd (f)^2 = \int_{\D} |f'(z)|^2 \, dA(z) < \infty,
            \]
            where $dA(z)$ represents the normalized area measure of $\D$. Actually, $\cd$ is the unique Hilbert space with this property \cite{AF} and it is easy to see that the inclusion $\cd \subset \cb$ is proper. In the case that $f$ is univalent, the change of variables $\zeta = f(z)$ shows that $\rho_{\cd}^2(f)$ is precisely the normalized area of $f(\D)$. If $f(z) = \sum_{k \geq 0} a_k z^k$, Parseval's identity yields the alternative representation
            \[
                \rho_\cd^2 (f) = \sum_{k = 1}^\infty k |a_k|^2.
            \]
        
        We refer the reader to the seminal paper \cite{ACP} and the monograph \cite{EKMR} for further information about $\cb$ and $\cd$, respectively.

        In 1978, Wirths \cite{W} showed that if $n \geq 2$ and $f(z) = \sum_{k \geq 0} a_k z^k$ is a normalized Bloch function (\emph{i.e.}, a function with $\| f \|_{\cb} = 1$), then
            \begin{equation} \label{eqn:Wirths}
                |a_n| \leq \dfrac{n+1}{2n} \left( \dfrac{n+1}{n-1} \right)^\frac{n-1}{2} =: B_n.
            \end{equation}
            Furthermore, equality is possible if and only if $f$ is a rotation of $z \mapsto B_n z^n$ (see \cite[Corollary~2]{KW3}). If $n = 1$, it trivially follows from \eqref{eqn:BlochSeminorm} that $|a_1| \leq 1$ so we define $B_1 = 1$ accordingly.

        Related to this bound, consider the quadratic functional
            \[
                \Phi_n(f) = \sum_{k = 1}^n k |a_k|^2, \quad f \in \cb.
            \]
            In other words, $\Phi_n (f)$ is the truncation of the (perhaps divergent) series defining $\rho_\cd^2 (f)$. In 2017, Kayumov \cite{Kayumov} proved that
            \begin{equation} \label{eqn:KayumovWirths}
                \Phi_n (f) \leq nB_n^2
            \end{equation}
            for all $f \in \cb$ with $\| f \|_\cb \leq 1$ when $n \leq 5$. Actually, his proof yields that
            \[
                \sum_{k = [n/2]}^n k |a_k|^2 \leq nB_n^2,
            \]
            where $[x]$ denotes the integer part of $x$, for all $n$ \cite[p.~467]{Kayumov}.
            These results were later complemented by Kayumov and Wirths \cite{KW1}, who showed that \eqref{eqn:KayumovWirths} holds for $n = 6$ as well. Other related extremal problems for coefficients of Bloch functions have been studied since then (e.g., \cite{KW2} and \cite{Khasianov}).

        In \cite{EMV}, it was proved that \eqref{eqn:KayumovWirths} holds for all $n$ when the coefficients of $f$ are non-negative. Furthermore,
            \[
                n B_n^2 \leq \max \left \{ \Phi_n(f) : \| f\|_\cb = 1 \right\} \leq \dfrac{n^n}{(n-1)^{n-1}} < \dfrac{4}{e} nB_n^2,
            \]
            so the bound $n B_n^2$ is of the correct order of magnitude. Even though it was not made explicit in the original paper \cite{Kayumov}, it was mentioned in \cite[p.~6]{KW4} that \eqref{eqn:KayumovWirths} could be valid for every $n$.

        Despite this evidence, the first main result of this work is to show that \eqref{eqn:KayumovWirths} does not hold when $n$ is large enough.

        \begin{thm} \label{thm:FailureLargeN}
            Let $n \geq 27$. Then, there exists $f \in \cb$ with $\| f \|_\cb = 1$ such that
                \[
                    \Phi_n (f) > n B_n^2.
                \]
        \end{thm}

        Nevertheless, we are able to show that \eqref{eqn:KayumovWirths} holds for all $n \leq 9$. Namely, we use a refinement of the methods of \cite{Kayumov} to prove the following:

        \begin{thm} \label{thm:KWforNsmall}
            Suppose that $6 \leq n \leq 9$. Then we have
                \[
                    \Phi_n (f) \leq n B_n^2
                \]
                for all normalized Bloch functions. Furthermore, equality holds if and only if $f(z) = \lambda B_n z^n$ for some unimodular constant $\lambda$.
        \end{thm}

        Observe that, as stated before, the case $n = 6$ was already known \cite{KW1}. However, the present argument will represent a new proof of this fact.

        To the best of our knowledge, the cases $10 \leq n \leq 26$ are still open.

        This work is organized as follows: the construction of the counterexamples to \eqref{eqn:KayumovWirths} which lead to the proof of Theorem~\ref{thm:FailureLargeN} can be found in Section~\ref{sec:FailureOfTheConjecture}. Indeed, we will present one function for the case $n = 27$ and another that covers all other cases $n \geq 28$.  Section~\ref{sec:PositiveSolutions} is devoted to the positive solution for $n \in \{6, 7, 8, 9\}$. In Section~\ref{sec:ConcludingRemarks}, we extend Theorem~\ref{thm:FailureLargeN} to the family of generalized weighted Dirichlet functionals $\Phi_n^t$. Finally, an appendix with a series of elementary but technical results is included for the correct understanding of the main sections.

        \subsection*{AI statement} This work has been carried out with assistance from Claude Opus 5.5 (Anthropic). It began with the aim of continuing the work \cite{EMV} and attempting to prove \eqref{eqn:KayumovWirths} for every $n$. When asked about the case $n = 7$, it produced a first draft of solution for $6 \leq n \leq 9$. In particular, it proposed the use of Lemma~\ref{lemma2} and the restricted parameter $\kappa$ that can be found in Section~\ref{sec:PositiveSolutions}. While exploring a suitable replacement of Lemma~\ref{lemma2} for larger $n$, it found a counterexample for \eqref{eqn:KayumovWirths} when $n \geq 28$. Further discussions led to the failure of the case $n = 27$ as well.
        All these proofs have been distilled, simplified and written in their final form by the authors.
    
        \subsection*{Acknowledgements}
    
            The three authors were partially supported by the grant PID2025-168349NB-I00 from Ministerio de Ciencia, Innovación y Universidades and by the project CIGE/2026/195 from Conselleria d'Educació, Cultura i Universitats, Generalitat Valenciana. The work of the second author was also supported by the grant PID2024-160326NA-I00 from the same ministry.

    \section{Failure of \eqref{eqn:KayumovWirths} for large $n$} \label{sec:FailureOfTheConjecture}

        In this section we will construct the counterexample needed to prove Theorem~\ref{thm:FailureLargeN}. In light of \cite[Theorem~2]{EMV}, the coefficients of this function must not have constant argument. Actually, the same proof yields that it is not possible that
            \[
                \arg a_k = \alpha + k \beta, \quad 1 \leq k,
            \]
            for certain $\alpha, \beta \in \mathbb{R}$. This represents the main difficulty of the following lemma.

        \begin{lem}\label{lem:NormCounterexample}
            Let $n \geq 4$. Consider the function
                \begin{equation} \label{eqn:Counterexample}
                    f_\varepsilon (z) = \varepsilon z + B_n z^n - \dfrac{\varepsilon}{2n-1} \left( \dfrac{2n}{n+1} \right)^{2} B_n^2 z^{2n-1}, \quad \varepsilon > 0.
                \end{equation}
                Then, there exists $\varepsilon_n > 0$ such that $\rho_\cb^2 (f_\varepsilon) = 1 + \dfrac{16 \varepsilon^2}{(n+1)^2}$ for all $0 < \varepsilon < \varepsilon_n$.
        \end{lem}

        \begin{proof}
            Write $z = r e^{it}$. Then we have that
                \begin{eqnarray}
                    |f_\varepsilon'(z)|^2 & = & \left| nB_n z^{n-1} + \varepsilon \left(1 - \left(\dfrac{2n}{n+1}\right)^2 B_n^2 z^{2n-2} \right)\right|^2 \nonumber \\
                        & = & n^2 B_n^2 r^{2n-2} + 2 nB_n \varepsilon\operatorname{Re} \left\{ z^{n-1} \left( 1 - \left(\dfrac{2n}{n+1}\right)^2 B_n^2 \overline{z}^{2n-2} \right) \right\} + \varepsilon^2 \left| 1 - \left(\dfrac{2n}{n+1}\right)^2 B_n^2 z^{2n-2} \right|^2 \nonumber \\
                        & = & n^2 B_n^2 r^{2n-2} + 2n B_n A(n, r, t) \varepsilon + B(n, r, t) \varepsilon^2, \label{eqn:aux1}
                \end{eqnarray}
                where
                \begin{eqnarray}
                    A(n, r, t) & = & \left( 1 - \left( \dfrac{2n}{n+1}\right)^2 B_n^2 r^{2n-2} \right) r^{n-1} \cos \big( (n-1)t \big) \nonumber \\
                    & = & \left(1 - \left( \dfrac{n+1}{n-1} r^2\right)^{n-1}\right) r^{n-1} \cos \big( (n-1)t \big) ,\label{eqn:auxA} \\
                    B(n, r, t) & = & 1 - 2 \left(\dfrac{2n}{n+1} \right)^2 B_n^2 r^{2n-2} \cos \big( 2(n-1)t\big) + \left(\dfrac{2n}{n+1} \right)^4 B_n^4 r^{4n-4} \nonumber \\
                    & = & \left( 1 + \left( \dfrac{n+1}{n-1} r^2\right)^{n-1} \right)^2 - 4 \left( \dfrac{n+1}{n-1} r^2\right)^{n-1} \cos^2 \big((n-1)t\big) \label{eqn:auxB}
                \end{eqnarray}
                (recall the definition of $B_n$ \eqref{eqn:Wirths}). If $r > 0$ and $\varepsilon$ are fixed, observe from \eqref{eqn:aux1}, \eqref{eqn:auxA} and \eqref{eqn:auxB} that $|f_{\varepsilon}'(re^{it})|^2$ is a concave quadratic function of $\cos \big( (n - 1)t \big)$. That quadratic function attains its maximum at its vertex
                \[
                    v_\varepsilon (r) = \dfrac{n+1}{8 \varepsilon} \left( \left( \frac{n+1}{n-1} r^2\right)^{-\frac{n-1}{2}}-\left( \frac{n+1}{n-1} r^2\right)^\frac{n-1}{2}\right).
                \]
                For the critical radius $r_n = \sqrt{\frac{n-1}{n+1}}$ we have that $v_\varepsilon(r_n) = 0$ and therefore on that radius we deduce that
                \begin{eqnarray}
                    |f_{\varepsilon}'(r_n e^{it})|^2 (1 - r_n^2)^2 & \leq & \left( n^2 B_n^2 r_n^{2n-2} + \left( 1 + \left( \dfrac{n+1}{n-1} r_n^2 \right)^{n-1}\right)^2 \varepsilon^2 \right) (1 - r_n^2)^2 \label{eqn:aux2} \\
                    & = & 1 + \dfrac{16 \varepsilon^2}{(n+1)^2} \nonumber,
                \end{eqnarray}
                and equality in \eqref{eqn:aux2} is possible when $e^{(n-1)ti}$ is a purely imaginary number. This in particular shows that $\rho_\cb^2 (f_\varepsilon) \geq 1 + \frac{16 \varepsilon^2}{(n+1)^2}$ for every $\varepsilon > 0$.

            We will now show that the converse inequality holds if $\varepsilon$ is small enough. If $r = 0$, it is clear that $|f_{\varepsilon}'(0)|^2 = \varepsilon^2 < 1 + \frac{16 \varepsilon^2}{(n+1)^2}$ if $\varepsilon < 1$. From now on assume that $0 < r < 1$ and thus $v_\varepsilon(r)$ is well defined. In order to simplify the later notation, write $R = \left( \frac{n+1}{n-1} r^2\right)^{n-1}$ so that
                \[
                    v_\varepsilon(r) = \dfrac{n+1}{8 \varepsilon} \left( \sqrt{\dfrac{1}{R}}-\sqrt{R} \right).
                \]
                Observe that $v_\varepsilon(r)$ is a continuous and strictly decreasing function of $R$ (and therefore of $r$) and monotonic with respect to $\varepsilon$.
                
            Even if there is no $t$ such that $\cos \big( (n - 1) t \big) = v_\varepsilon(r)$, the following inequality holds:
                \begin{eqnarray*}
                    |f_{\varepsilon}'(re^{it})|^2 & \leq & \dfrac{(n+1)^2}{4}  R + \dfrac{(n+1)^2 (1 - R)^2}{8} - \dfrac{(n+1)^2}{16} R \left( \sqrt{R} - \sqrt{\dfrac{1}{R}}\right)^2 + (1 + R)^2 \varepsilon^2 \nonumber \\
                        & = & \dfrac{(n+1)^2 (1 + R)^2}{16} + (1 + R)^2 \varepsilon^2 = \left( 1 + \dfrac{16 \varepsilon^2}{(n+1)^2} \right) \dfrac{(n + 1)^2 (1 + R)^2}{16}. 
                \end{eqnarray*}
                
            According to Remark~\ref{rem:AppendixRemark1} in the Appendix, there exists $r_0 < r_n$ independent of $\varepsilon$ such that $(1 + R)(1 - r^2) \leq \frac{4}{n+1}$   if $r_0 \leq r < 1$. We split the proof by cases.

            \begin{enumerate}
                \item If $r \geq r_0$, then
                    \[
                        |f_\varepsilon'(re^{it})|^2 (1 - r^2)^2 \leq \left( 1 + \dfrac{16 \varepsilon^2}{(n+1)^2} \right) \dfrac{(n + 1)^2 (1 + R)^2 (1 - r^2)^2}{16} \leq 1 + \dfrac{16 \varepsilon^2}{(n+1)^2}.
                    \]
                    
                \item When $0 < r < r_0$, in particular $v_\varepsilon(r) > v_\varepsilon(r_0) > v_\varepsilon(r_n) = 0$. Pick $\varepsilon_n < 1$ such that $v_{\varepsilon_n} (r_0) \geq 1$. For every $0 < \varepsilon < \varepsilon_n$ the quantity $|f_{\varepsilon}'(re^{it})|^2$ is then maximized when $e^{(n-1)ti} = 1$. Thus
                    \[
                        |f_{\varepsilon}'(re^{it})|^2 \leq \dfrac{(n+1)^2}{4} R + (n+1) \sqrt{R} (1 - R) \varepsilon + \big( (1 + R)^2 - 4R \big) \varepsilon^2 = \left( \dfrac{n+1}{2} \sqrt{R} + (1 - R) \varepsilon \right)^2.
                    \]
                    Summing up, we have that 
                    \[
                        |f_{\varepsilon}'(re^{it})| (1 - r^2) \leq \left( \dfrac{n+1}{2} \sqrt{R} + (1 - R) \varepsilon \right) (1 - r^2) \leq nB_n r^{n-1} (1 - r^2) + \varepsilon,
                    \]
                    since $R < 1$ because $r < r_0 < r_n$. The monotonicity of the function $r^{n-1} (1 - r^2)$ (see Lemma~\ref{lem:AppendixLemma2} in the Appendix) yields that we may also choose $\varepsilon_n$ such that $nB_n r_0^{n-1} (1 - r_0^2) + \varepsilon_n < 1$. Under this new assumption, we have that $|f_{\varepsilon}'(re^{it})| (1 - r^2) < 1$ if $0 < r < r_0$. This concludes the proof.
            \end{enumerate}
        \end{proof}

        \begin{proof}[Proof of Theorem~\ref{thm:FailureLargeN} for $n \geq 28$]
            Let $f_\varepsilon$ be the function defined in \eqref{eqn:Counterexample} and assume that $0 < \varepsilon < \varepsilon_n$. Then we have that
                \begin{eqnarray*}
                    \Phi_n \left( \dfrac{f_\varepsilon}{\| f_\varepsilon \|_\cb} \right) & = & \dfrac{\varepsilon^2 + n B_n^2}{1 + \frac{16 \varepsilon^2}{(n+1)^2}} \\
                        & = & n B_n^2 \dfrac{1 +\frac{4n \varepsilon^2}{(n+1)^2} \left(\frac{n-1}{n+1}\right)^{n-1}}{1 + \frac{16 \varepsilon^2}{(n+1)^2}}.
                \end{eqnarray*}
                Thus, it is enough to check that
                \begin{equation} \label{eqn:aux4}
                    n \left( \dfrac{n-1}{n+1}\right)^{n-1} > 4.
                \end{equation}
                The proof concludes by observing that the sequence in the left-hand side of \eqref{eqn:aux4} is increasing and that $28 \left( \frac{27}{29} \right)^{27} > 4$.
        \end{proof}

        Observe that $n=28$ is the smallest integer such that \eqref{eqn:aux4} holds. Hence, a counterexample for $n = 27$ must be necessarily different. To this end, we will add an extra high-order term which will not contribute to $\Phi_n$.

        \begin{proof}[Proof of Theorem~\ref{thm:FailureLargeN} for $n = 27$]
            Given $\varepsilon > 0$, consider the polynomial
                \begin{equation} \label{eqn:Countexample27}
                    f_\varepsilon (z) = 14 \varepsilon z + B_{27} z^{27} - \dfrac{14}{53} \left( \dfrac{14}{13}\right)^{26} \varepsilon z^{53} + \dfrac{1}{79} \left( \dfrac{14}{13} \right)^{39} \varepsilon^2 z^{79}.
                \end{equation}
                Thus,
                \[
                    \Phi_{27} \left( \dfrac{f_\varepsilon}{\| f_\varepsilon \|_\cb}\right) = 27 B^2_{27} \dfrac{1 + 27 \left( \frac{13}{14} \right)^{26} \varepsilon^2}{\| f_\varepsilon \|^2_\cb}.
                \]
                We will check that $\| f_\varepsilon \|^2_\cb < 1 + 27 \left( \frac{13}{14} \right)^{26} \varepsilon^2$ if $\varepsilon$ is small enough.

            It is plain that
                \[
                    f_\varepsilon' (z) - 27 B_{27} z^{26} = 14 \varepsilon - 14 \left( \dfrac{14}{13} \right)^{26} \varepsilon z^{52} + \left( \dfrac{14}{13} \right)^{39} \varepsilon^2 z^{78}.
                \]
                Therefore, $f_\varepsilon \to B_{27} z^{27}$ in the Bloch topology when $\varepsilon \to 0^+$. In particular, if
                \[
                    |f_\varepsilon' (z_\varepsilon)| (1 - |z_\varepsilon|^2) = \rho_\cb (f_\varepsilon)
                \]
                then $|z_\varepsilon|^2 \to \frac{13}{14}$ (see Lemma~\ref{lem:AppendixLemma2}). In particular, if $\varepsilon$ is small enough we may use Remark~\ref{rem:AppendixRemark2} to see that
                \[
                    \|f_\varepsilon \|_\cb^2 = \sup_{r_1 \leq |z| \leq r_2} \left \{ |f_\varepsilon' (z)|^2 (1 - |z|^2)^2 \right \} < 1 + 27 \left( \dfrac{13}{14}\right)^{26} \varepsilon^2,
                \]
                so $f_\varepsilon$ is a valid counterexample of \eqref{eqn:KayumovWirths}.
        \end{proof}

    \section{Positive solutions for $6 \leq n \leq 9$} \label{sec:PositiveSolutions}

        If $\| f \|_\cb \leq 1$, the definition of the Bloch seminorm yields that
            \[
                |f'(z)|^2 \leq \dfrac{1}{(1 - r^2)^2}, \quad |z| = r.
            \]
            Integrating over the circle $\{ |z| = r \}$ and applying Parseval's identity imply that
            \begin{equation} \label{eqn:Parseval}
                \sum_{k = 1}^\infty k^2 |a_k|^2 r^{2k-2} = \dfrac{1}{2\pi} \int_0^{2\pi} |f'(re^{it})|^2 \, dt \leq \dfrac{1}{(1 - r^2)^2}, \quad 0 < r < 1.
            \end{equation}
            This elementary but ingenious inequality allows us to deduce \eqref{eqn:KayumovWirths} for $1 \leq n \leq 5$ when substituting by the critical radius $r_n=\sqrt{\frac{n-1}{n+1}}$ (see \cite{Kayumov}). However, it is not sufficient by itself to cover the cases $6 \leq n \leq 9$.

        We will divide the analysis into two complementary propositions based on the relative magnitude of the $n$-th Taylor coefficient $|a_n|$ with respect to the first one $|a_1|$.

        We first consider the case where the $n$-th coefficient is relatively small. In this case, the defect generated by the intermediate coefficients combined with Parseval's identity is sufficient to establish the desired inequality under some additional technical hypotheses which will be checked later in Proposition~\ref{prop:n=6,7,8,9}.
        
        \begin{prop}\label{prop:afirmative1}
            Let $n \geq 3$. Consider the sequence
                \begin{equation} \label{eqn:Weights}
                    w_k = k \left( \dfrac{n-1}{n+1}\right)^{k-1}, \quad 1 \leq k, \tag{$w_k$}
                \end{equation}and define
                \begin{equation} \label{eqn:WMin}
                    w_{\min} = \min_{2 \leq k \leq n-1} \{ w_k \}. \tag{$w_{\min}$}
                \end{equation} Assume that there exists a number $\kappa$ such that the conditions
                \begin{equation}
                    \tag{H0} 0 < \kappa < \dfrac{1}{w_n} \quad \mbox{and} \quad w_{\min} > w_n > 1;\label{eqn:H0}
                \end{equation}
                and
                \begin{equation}
                    \tag{H1} \dfrac{w_{\min} - 1}{w_{\min} - w_n}+ \dfrac{4 (1 - w_n)^2}{\kappa^2 w_n^3}\leq n B_n^2. \label{eqn:H3}
                \end{equation}
                are satisfied.

             If $f(z)=\sum_{k \ge 0} a_k z^k$ with $\norm{f}_{\cb}\le 1$ satisfies
                \begin{equation} \label{eqn:SmallAN}
                    n |a_n|^2 < \dfrac{4 (1 - w_n)^2}{\kappa^2 w_n^3} |a_1|^2,
                \end{equation}
                then
                \[
                    \sum_{k=1}^n k |a_k|^2 < n B_n^2.
                \]
        \end{prop}

        \begin{proof}
            In order to simplify the notation, write
                \begin{equation} \label{eqn:SequenceYk}
                    y_k = \sqrt{kw_k} |a_k|, \quad k \geq 1,
                \end{equation}
                so that
                \[
                    \sum_{k=1}^n k |a_k|^2 = \sum_{k=1}^{n} \frac{y_k^2}{w_k}
                \]
                By the definition of $w_{\min}$, we will have that
                \begin{equation} \label{eqn:aux9}
                    \sum_{k=1}^n k |a_k|^2 \le y_1^2 + \frac{y_n^2}{w_n} + \frac{1}{w_{\min}} \sum_{k=2}^{n-1} y_k^2.
                \end{equation}
                
            Applying \eqref{eqn:Parseval} at $r_n = \sqrt{\frac{n-1}{n+1}}$ we see that
                \[
                    \sum_{k = 2}^{n-1} y_k^2 \le n w_n B_n^2 - y_1^2 - y_n^2.
                \]
                Then, substituting in \eqref{eqn:aux9} gives
                \[
                    \sum_{k=1}^n k |a_k|^2 \le \frac{w_n}{w_{\min}} nB_n^2 + \frac{w_{\min} - 1}{w_{\min}} y_1^2 + \frac{w_{\min} - w_n}{w_n w_{\min}} y_n^2.
                \]  
                Factoring out $\frac{w_{\min} - w_n}{w_{\min}}$ and using $y_1 \le 1$ together with \eqref{eqn:SmallAN} and \eqref{eqn:H0}, we obtain the strict inequality
                \[
                    \sum_{k=1}^n k |a_k|^2 < \frac{w_n}{w_{\min}} n B_n^2 + \frac{w_{\min} - w_n}{w_{\min}} \left( \frac{w_{\min} - 1}{w_{\min} - w_n} + \dfrac{4 (1 - w_n)^2}{\kappa^2 w_n^3} \right).
                \]
                Applying hypothesis \eqref{eqn:H3} guarantees that $\sum_{k=1}^n k |a_k|^2 < n B_n^2$.
        \end{proof}

        We now turn to the complementary statement of Proposition~\ref{prop:afirmative1} where the $n$-th coefficient is dominant. To this end, we first need the following refinement of \eqref{eqn:Parseval}.

        \begin{lem}\label{lemma2}
            Let $r \in (0, 1)$ and $m \ge 1$ be an integer. If $f(z)=\sum_{k \ge 0} a_k z^k$ and $\norm{f}_{\cb}\le 1$, then
            \[
                \sum_{k = 1}^\infty k^2 |a_k|^2 r^{2k-2} + \left| \sum_{j = 1}^\infty j(j+m) a_{j+m} \overline{a_j} r^{2j+m-2} \right| \le \frac{1}{(1-r^2)^2}.
            \]
        \end{lem}

        \begin{proof}
            Consider $0<r<1$ and the real-valued function $G(t) = |f'(re^{it})|^2$, which is continuous on $[0, 2\pi]$ and clearly satisfies
            \[
                0 \le G(t) \le \frac{1}{(1-r^2)^2}, \quad \forall t \in [0, 2\pi].
            \]
            
            Writing $G(t) = f'(re^{it}) \overline{f'(re^{it})}$ and taking into account that the power series of $f'(z)$ converges absolutely on $\{|z| = r \}$, we obtain that the Fourier coefficients $(c_l)_{l \in \Z}$ of $G$ are given by
            \[
                c_0 = \sum_{k=1}^\infty k^2 |a_k|^2 r^{2k-2} \quad \text{and} \quad c_l = \overline{c_{-l}} = \sum_{j=1}^\infty j(j+l) a_{j+l} \overline{a_j} r^{2j+l-2}, \quad l \ge 1.
            \]

            Now, fix $\varphi \in \mathbb{R}$ and consider the trigonometric polynomial $P_\varphi(t) = 1 + \cos(mt - \varphi)$. Observe that $P_\varphi$ is non-negative and, furthermore,
                \[
                    \dfrac{1}{2\pi} \int_{0}^{2\pi} P_\varphi (t) \, dt = 1.
                \]
            On the one hand, we get
                \begin{equation} \label{eqn:aux8}
                    \frac{1}{2\pi} \int_0^{2\pi} G(t) \big(1 + \cos(mt - \varphi)\big) \, dt \le \frac{1}{(1-r^2)^2}.
                \end{equation}
                On the other hand, taking into account that
                \[
                    \frac{1}{2\pi} \int_0^{2\pi} G(t) \cos(mt - \varphi) dt= \mathrm{Re}\left( e^{i\varphi} \frac{1}{2\pi} \int_0^{2\pi} G(t) e^{-imt} \, dt \right) = \mathrm{Re}\left( c_m e^{i\varphi} \right),
                \]
                expressing the integral on the left-hand side of \eqref{eqn:aux8} via Fourier coefficients implies that
                \[
                    c_0 + \mathrm{Re}\left( c_m e^{i\varphi} \right) \le \frac{1}{(1-r^2)^2}.
                \]
            
            Choosing  $\varphi = -\arg(c_m)$ we get  $\mathrm{Re}\left( c_m e^{i\varphi} \right) = |c_m|$. Substituting the explicit series representations of $c_0$ and $c_m$ yields
                \[
                    \sum_{k=1}^\infty k^2 |a_k|^2 r^{2k-2} + \left| \sum_{j=1}^\infty j(j+m) a_{j+m} \overline{a_j} r^{2j+m-2} \right| \le \frac{1}{(1-r^2)^2},
                \]
                which completes the proof.
        \end{proof}
        \begin{lem}\label{lem:applyImproveParseval}
            Let $n \geq 2$. Consider the sequence $(w_k)_k$ defined in \eqref{eqn:Weights}. For every $f (z) = \sum_{k \geq 0} a_k z^k$ in $\mathcal{B}$ with $\|f\|_{\mathcal{B}} \le 1$, we have
                \[
                    w_n \left( n B_n^2  - \Phi_n(f)\right) \ge \sum_{k=2}^{n-1} k(w_k-w_n)|a_k|^2 + \sum_{k=n+1}^\infty k w_k  |a_k|^2 + \left|R\right| - (w_n - 1) |a_1|^2,
                \]
                where $R=\sum_{j \geq 1} j(j+n-1) a_{j+n-1} \overline{a_j} r_n^{2j+n-3}$.
        \end{lem}

        \begin{proof}
            Applying Lemma \ref{lemma2} for $m = n - 1$ at $r_n = \sqrt{\frac{n-1}{n+1}}$, we have
            \[
                \sum_{k=1}^\infty k w_k |a_k|^2 + |R| \le \frac{1}{(1-r_n^2)^2} = n w_n B_n^2.
            \]
            Subtracting $w_n \Phi_n(f)$ from both sides yields the desired inequality since $w_1=1$.
        \end{proof}

        Similarly to Proposition~\ref{prop:afirmative1}, we will need to control the interaction between $a_1$, the dominant $a_n$ and $a_{2n-1}$. Observe that this last coefficient does not play a role for $\Phi_n(f)$ but it is present in the remainder $R$ introduced in Lemma~\ref{lem:applyImproveParseval}. This task will require delicate estimates based on further technical hypotheses which we can safely assume when $6 \leq n \leq 9$.

        \begin{prop}\label{prop:afirmative2}
            Let $n \geq 3$ and consider the sequence introduced in \eqref{eqn:Weights} and the quantity $w_{\min}$ defined in \eqref{eqn:WMin}. Define
                \begin{equation}
                    w_{\max} = \max_{2 \leq k \leq n-1} \{ w_k \}. \tag{$w_{\max}$}
                \end{equation}
                Assume that there exists a number $\kappa$ satisfying \eqref{eqn:H0} and the additional conditions
                \begin{equation} \label{eqn:H1}
                    \dfrac{\kappa w_n}{2} + \dfrac{\kappa^2 w_n^2 w_{\max}}{4(w_{\min} - w_n)} < 1 \tag{H2}
                \end{equation}
                and
                \begin{equation} \label{eqn:H2}
                    \dfrac{\kappa w_n}{2} + w_n - 1 < 1. \tag{H3}
                \end{equation}
            
            Suppose that $f(z)=\sum_{k \ge 0} a_k z^k$ with $\norm{f}_{\cb}\le 1$ satisfies
                \begin{equation}\label{condition_lemma}
                    n |a_n|^2 \geq \dfrac{4 (1 - w_n)^2}{\kappa^2 w_n^3} |a_1|^2,
                \end{equation}
                then
                    \[
                        \sum_{k=1}^n k |a_k|^2 \leq n B_n^2.
                    \]
                 Furthermore, equality holds if and only if $f(z) = \lambda B_n z^n$ for a certain unimodular constant $\lambda$.   
        \end{prop}

        \begin{proof}
            Recall Lemma \ref{lem:applyImproveParseval}:
            \begin{equation*} 
                w_n \left( n B_n^2  - \Phi_n(f) \right) \ge \sum_{k=2}^{n-1} k(w_k-w_n)|a_k|^2 + \sum_{k=n+1}^\infty k w_k |a_k|^2 + |R| - (w_n - 1)|a_1|^2.
            \end{equation*}
            If we can prove that the right-hand side of this inequality is non-negative, we are done.

            Using the notation \eqref{eqn:SequenceYk} and triangle inequality we get 
                \[
                    |R| \geq y_1y_n-y_ny_{2n-1}- S,
                \]
                where $S=\sum_{k=2}^{n-1} y_ky_{k+n-1}+\sum_{k\geq n+1} y_ky_{k+n-1}$.
            
            On the other hand, we will see that
            \[
                \sum_{k=2}^{n-1} k(w_k-w_n)|a_k|^2 + \sum_{k=n+1}^\infty k w_k |a_k|^2 \geq \kappa w_n S+(w_n-1)y_{2n-1}^2
            \]
            which, in terms of the notation \eqref{eqn:SequenceYk}, is 
                \begin{equation}\label{eq_domitaion}
                    \sum_{k=2}^{n-1} \frac{w_k - w_n}{w_k} y_k^2 + \sum_{k=n+1}^\infty y_k^2 \geq \kappa w_n S+(w_n-1)y_{2n-1}^2
                \end{equation}
                We will prove this completing squares and using the hypotheses. We split the computations in four steps depending on the range of $k$.
                
            Firstly, observe that 
                \begin{align}
                    \sum_{k=2}^{n-1} \frac{w_k - w_n}{w_k} y_k^2 -&\kappa w_n \sum_{k=2}^{n-1} y_ky_{k+n-1}\nonumber\\
                    &= \sum_{k=2}^{n-1}\left( \sqrt{\dfrac{w_k - w_n}{w_k}}y_k - \frac{\kappa w_n}{2} \sqrt{\dfrac{w_k}{w_k-w_n}}y_{k+n-1}\right)^2- \sum_{k=2}^{n-1}\frac{\kappa^2 w_n^2 w_k}{4(w_k - w_n)}y_{k+n-1}^2\label{first_step}\\
                    &\geq - \frac{\kappa^2 w_n^2 w_{\max}}{4(w_{\min} - w_n)} \sum_{k=2}^{n-1}y_{k+n-1}^2\nonumber.
                \end{align}
                
            Secondly, using hypothesis \eqref{eqn:H1} and the change of index $k = j + n - 1$ we see that
                \begin{equation}\label{second_step}
                    \sum_{k=n+1}^{2n-2}y_k^2\geq \left( \frac{\kappa w_n}{2}+\frac{\kappa^2 w_n^2 w_{\max}}{4(w_{\min} - w_n)}\right)\sum_{k=n+1}^{2n-2}y_k^2 = \frac{\kappa w_n}{2}\sum_{k=n+1}^{2n-2}y_k^2 + \frac{\kappa^2 w_n^2 w_{\max}}{4(w_{\min} - w_n)}\sum_{j=2}^{n-1}y_{j+n-1}^2.
                 \end{equation}
                
            Thirdly, using \eqref{eqn:H2} we get
                \begin{equation}\label{third_step}
                    y_{2n-1}^2 \geq (w_n-1) y_{2n-1}^2+ \dfrac{\kappa w_n}{2} y_{2n-1}^2
                \end{equation}
                and \eqref{eqn:H0} and another change of index yield
                \begin{equation}\label{fourth_step}
                    \sum_{k = 2n}^\infty y_k^2\geq \kappa w_n \sum_{k = 2n}^\infty y_k^2=\dfrac{\kappa w_n}{2}\sum_{k = 2n}^\infty y_k^2 +\dfrac{\kappa w_n}{2}\sum_{j = n+1}^\infty y_{j+n-1}^2. 
                \end{equation}
                
            Finally, by combining all these inequalities, we obtain  
                \begin{align*}
                    \sum_{k=2}^{n-1} \frac{w_k -w_n}{w_k} y_k^2 + \sum_{k=n+1}^\infty y_k^2 & \geq  \kappa w_n \sum_{k=2}^{n-1} y_ky_{k+n-1}+  (w_n-1) y_{2n-1}^2 + \frac{\kappa w_n}{2}\sum_{k = n+1}^\infty (y_k^2+y_{k+n-1}^2)\\
                &\geq  \kappa w_n \sum_{k=2}^{n-1} y_ky_{k+n-1}+  (w_n-1) y_{2n-1}^2 + \kappa w_n\sum_{k= n+1}^\infty y_ky_{k+n-1}
                \end{align*}
                that yields \eqref{eq_domitaion}. Therefore, using the non-negativity of $|R|$, $y_1=|a_1|$ and Lemma \ref{lem:auxiliary_max} in the Appendix, we obtain that
                \begin{align*}
                    w_n \left( n B_n^2  - \Phi_n (f) \right) &\geq \max \left\{ 0, y_1y_n-y_ny_{2n-1}- S\right\} +\kappa w_n S+(w_n-1)(y_{2n-1}^2-y_1^2)\\
                    & \geq \kappa w_n y_n \max\left\{ 0, y_1- y_{2n-1}\right\}+(w_n-1)(y_{2n-1}^2-y_1^2).
                \end{align*}
                
            Hence, if $y_{2n-1} \geq y_1$, 
                \[
                    w_n \left( n B_n^2  - \Phi_n (f)\right) \geq (w_n-1)(y_{2n-1}^2-y_1^2)\geq 0,
                \]
                and if $y_{2n-1}< y_1$ we have instead that
                \begin{align*}
                    w_n \left( n B_n^2  - \Phi_n (f) \right) &\geq \kappa w_n y_n(y_1-y_{2n-1}) +(w_n-1)(y_{2n-1}^2-y_1^2)\\
                    & =(\kappa w_n y_n + (1-w_n)(y_{2n-1}+y_1))(y_1-y_{2n-1})\\
                    &> (\kappa w_n y_n + 2(1-w_n)y_1)(y_1-y_{2n-1})\geq 0.
                \end{align*}
                In the last inequality we have used \eqref{condition_lemma}. 
                
                It remains to characterize the equality case $\Phi_n(f) = n B_n^2$. Suppose that $\Phi_n(f) = n B_n^2$, so that $w_n \left( n B_n^2 - \Phi_n (f) \right) = 0$. Since $y_{2n-1} < y_1$ gives a strict inequality above, we only need to analyze the case $y_{2n-1} \geq y_1$. We split up this analysis into two parts
                \begin{itemize}
                    \item If $y_1>0$, since $y_{2n-1}\geq y_1>0$ and the hypothesis \eqref{eqn:H2} is strict, we get that the inequality in  \eqref{third_step} is actually strict. Therefore, $w_n \left( n B_n^2 - \Phi_n (f) \right) > 0$. 
                    \item If $y_1=0$, by the previous case we can suppose that $y_{2n-1}=0$. Now, since hypotheses \eqref{eqn:H0} and \eqref{eqn:H1} are strict, we obtain the inequalities in \eqref{fourth_step} and  \eqref{second_step} force $y_{k}=0$ for $k\geq n+1$. Then, analysing the squares in \eqref{first_step} we get $y_k=0$ for $2\leq k\leq n-1$. Therefore, $\Phi_n(f) = n B_n^2$ holds if and only if $f(z) = \lambda B_n z^n$ for a certain unimodular constant $\lambda$.  
                \end{itemize}
                This finishes the proof.
        \end{proof}

        \begin{prop}\label{prop:n=6,7,8,9}
            For $n \in \{ 6,7,8,9 \}$ there exists $\kappa = \kappa(n)$ such that the hypotheses \eqref{eqn:H0}, \eqref{eqn:H3}, \eqref{eqn:H1} and \eqref{eqn:H2} are satisfied.
        \end{prop}

        \begin{proof}
            To complete the proof, we summarize the required numerical values in the Tables \ref{tb:Weights} and \ref{tb:Hypotheses}.
                
            \begin{table}[H]
                \centering
                \small
                \setlength{\tabcolsep}{3pt} 
                \renewcommand{\arraystretch}{2.2} 
                
                \begin{tabular}{c !{\vrule width 1.2pt} c c c c c c c c c c c}
                \hline
                $n$ & $w_1$ & $w_2$ & $w_3$ & $w_4$ & $w_5$ & $w_6$ & $w_7$ & $w_8$ & $w_9$ & $w_{\max}$ & $w_{\min} - w_n$ \\
                \Xhline{1.2pt}
                $6$ & $1$ & $\dfrac{10}{7}$ & $\dfrac{75}{49}$ & $\dfrac{500}{343}$ & $\dfrac{3125}{2401}$ & $\dfrac{18750}{16807}$ & & & & $\dfrac{75}{49}$ & $\dfrac{3125}{16807}$ \\
                $7$ & $1$ & $\dfrac{3}{2}$ & $\dfrac{27}{16}$ & $\dfrac{27}{16}$ & $\dfrac{405}{256}$ & $\dfrac{729}{512}$ & $\dfrac{5103}{4096}$ & & & $\dfrac{27}{16}$ & $\dfrac{729}{4096}$ \\
                $8$ & $1$ & $\dfrac{14}{9}$ & $\dfrac{49}{27}$ & $\dfrac{1372}{729}$ & $\dfrac{12005}{6561}$ & $\dfrac{33614}{19683}$ & $\dfrac{823543}{531441}$ & $\dfrac{6588344}{4782969}$ & & $\dfrac{1372}{729}$ & $\dfrac{823543}{4782969}$ \\
                $9$ & $1$ & $\dfrac{8}{5}$ & $\dfrac{48}{25}$ & $\dfrac{256}{125}$ & $\dfrac{256}{125}$ & $\dfrac{6144}{3125}$ & $\dfrac{28672}{15625}$ & $\dfrac{131072}{78125}$ & $\dfrac{589824}{390625}$ & $\dfrac{256}{125}$ & $\dfrac{35176}{390625}$ \\[2mm]
                \hline
                \end{tabular}
                \caption{Values of $(w_j)_{1 \leq j \leq n}$ for $6 \leq n \leq 9$.}
                \label{tb:Weights}
            \end{table}
                
            \begin{table}[H]
                \small
                \setlength{\tabcolsep}{3pt} 
                \renewcommand{\arraystretch}{2.2} 
                \begin{tabular}{c !{\vrule width 1.2pt} c |c| c| c |c} \hline $n$ & $\kappa$ & $\frac{2(w_n - 1)}{\kappa w_n}$ & \eqref{eqn:H3} & \eqref{eqn:H1} & \eqref{eqn:H2} \\[1mm] \Xhline{1.2pt} $6$ & $\dfrac{1}{2}$ & $\dfrac{3886}{9375}$ & $1.77577<10.9806$ & $\dfrac{3028125}{3294172}<1$ & $\dfrac{13261}{33614}<1$ \\[2mm] $7$ & $\dfrac{1}{4}$ & $\dfrac{8056}{5103}$ & $4.38177<12.8426$ & $\dfrac{1617651}{4194304}<1$ & $\dfrac{13159}{32768}<1$ \\[2mm] $8$ & $\dfrac{1}{4}$ & $\dfrac{1805375}{823543}$ & $6.6811 < 14.7010$ & $\dfrac{1730263843}{3486784401}<1$ & $\dfrac{292102}{531441}<1$ \\[2mm] $9$ & $\dfrac{1}{4}$ & $\dfrac{199199}{73728}$ & $11.4974 < 16.5568$ & $\dfrac{214468927488}{214697265625}<1$ & $\dfrac{272927}{390625}<1$ \\[2mm] \hline \end{tabular}

                \caption{Verification of \eqref{eqn:H3}, \eqref{eqn:H1} and \eqref{eqn:H2} for $6 \leq n \leq 9$.}
                \label{tb:Hypotheses}
            \end{table}
        \end{proof}

         \begin{proof}[Proof of Theorem~\ref{thm:KWforNsmall}]

            The result follows directly by combining Propositions~\ref{prop:afirmative1}, \ref{prop:afirmative2} and \ref{prop:n=6,7,8,9}. 
        \end{proof}

        \begin{rk}
            The proof technique presented here reaches its natural limit at $n = 9$ and cannot be extended to $n \ge 10$. Indeed, hypothesis \eqref{eqn:H0} strictly requires that $w_{\min} > w_n$ and, for example, for $n=10$ we have
            \[
                w_{10} - w_2 = 10 \left(\frac{9}{11}\right)^9 - \frac{18}{11} \approx 0.0067.
            \]        
        \end{rk}

    \section{Generalized weighted Dirichlet functionals} \label{sec:ConcludingRemarks}

        In \cite{EMV}, the authors also investigated a generalized weighted truncated Dirichlet functional defined for any parameter $t \ge 0$ by
            \[
                \Phi_n^t(f) = \sum_{k=1}^n k^t |a_k|^2,
            \]
            and analyzed the validity of the corresponding bound
            \begin{equation}\label{eq:generalized_kayumov}
                \Phi_n^t(f) \le n^t B_n^2
            \end{equation}
            for all functions $f(z) = \sum_{k=0}^\infty a_k z^k$ in $\mathcal{B}$ satisfying $\|f\|_{\mathcal{B}} \le 1$. As noted in \cite{EMV}, if $t_1 \leq t_2$ then the validity of \eqref{eq:generalized_kayumov} for $t_1$ implies its validity for $t_2$. Hence, the inequality \eqref{eq:generalized_kayumov} holds unconditionally for all $t \geq 2$ by applying \eqref{eqn:Parseval}, while for functions with non-negative Taylor coefficients it remains true for $t \geq 1$ (see \cite{EMV}). However, when $t < 1$, the estimate fails for all sufficiently large $n$, establishing $t = 1$  as the critical threshold for this inequality among functions with non-negative coefficients.

        The function considered in Theorem \ref{thm:FailureLargeN} to disprove inequality \eqref{eqn:KayumovWirths} for large $n$ is also suitable as a counterexample for the generalized inequality \eqref{eq:generalized_kayumov} when $t < 2$.

        \begin{cor}\label{cor:FailureLargeNGeneralizedFuntional}
            Let $1 \leq t<2 $. Then, there exists $n_t \geq 2$ such that for every $n \geq n_t$ we have that
                \[
                    \Phi_n^t (f) > n^t B_n^2.
                \]
                for a function $f \in \cb$ with $\| f \|_\cb = 1$.
        \end{cor}

        \begin{proof}
            Let $n_t$ be such that
                \[
                    n^{2-t}\left(\dfrac{n-1}{n+1} \right)^{n-1} > 4, \quad n \geq n_t.
                \]
                The index $n_t$ certainly exists because the sequence is unbounded.
                
                Now, consider $n \geq n_t$ and $f_\varepsilon$ the function defined in \eqref{eqn:Counterexample} and assume that $0 < \varepsilon < \varepsilon_n$, where $\varepsilon_n$ is the constant defined in Lemma~\ref{lem:NormCounterexample}. Then, proceeding as in the proof of Theorem \ref{thm:FailureLargeN} we have that
                \[
                    \Phi_n^t \left( \dfrac{f_\varepsilon}{\| f_\varepsilon \|_\cb} \right)  = \dfrac{\varepsilon^2 + n^t B_n^2}{1 + \frac{16 \varepsilon^2}{(n+1)^2}}  =  n^t B_n^2 \dfrac{1 +\frac{4n^{2-t} \varepsilon^2}{(n+1)^2} \left(\frac{n-1}{n+1}\right)^{n-1}}{1 + \frac{16 \varepsilon^2}{(n+1)^2}}>n^t B_n^2,
                 \]
            which completes the proof.
        \end{proof}

    \appendix

        \section{Auxiliary lemmata}

            In this Appendix we include some technical results which are needed for the correct understanding of the work.

            \begin{lem} \label{lem:AppendixLemma2}
                Let $n \geq 2$. Then the function $g(x) = x^{\frac{n-1}{2}} (1 - x)$, $0 \leq x \leq 1$, is strictly increasing in $\left[0, \dfrac{n-1}{n+1} \right]$ and decreasing in $\left[\dfrac{n-1}{n+1}, 1 \right]$.
            \end{lem}

            \begin{proof}
                It is a straightforward consequence of the identity
                    \[
                        g'(x) = \dfrac{x^{\frac{n-3}{2}}}{2} \big( n - 1 - (n + 1) x\big).
                    \]
            \end{proof}

            \begin{lem} \label{lem:AppendixLemma1}
                Let $n \geq 4$. Consider the function
                    \[
                        g(x) = (1 + x^{n-1}) \left(1 - \frac{n-1}{n+1} x \right), \quad 0 \leq x \leq \dfrac{n+1}{n-1}.
                    \]
                    There exists $0 < x_0 < 1$ such that $g(x) \leq g(1) =\frac{4}{n+1}$ for all $x \in \left[x_0, \dfrac{n+1}{n-1}\right]$.
            \end{lem}

            \begin{rk} \label{rem:AppendixRemark1}
                Using the change of variables $r^2 = \frac{n-1}{n+1} x$, $r > 0$, the above lemma shows that there exists $r_0 < \sqrt{\frac{n-1}{n+1}} = r_n$ such that
                \[
                    \left( 1 + \left( \dfrac{n+1}{n-1} r^2 \right)^{n-1} \right) (1 - r^2) \leq\left( 1 + \left( \dfrac{n+1}{n-1} r_n^2 \right)^{n-1} \right) (1 - r_n^2) = \dfrac{4}{n+1}, \quad r_0 \leq  r \leq 1. 
                \]
            \end{rk}

            \begin{proof}
                A direct computation shows that
                    \[
                        g'(x) = \dfrac{n-1}{n+1} \big( (n+1) x^{n-2} - n x^{n-1} - 1 \big) \, \mbox{ and } \, g''(x) = \dfrac{n-1}{n+1} x^{n-3} \big( (n+1)(n-2) - n(n-1) x\big).
                    \]
                    Let $x_1 = \frac{(n+1)(n-2)}{n(n-1)}$ be the unique zero of $g''$. Observe that $x_1 < 1$ because of the fact that $(n+1)(n-2) = n^2 - n -2 < n(n-1)$.

                    Since $g'$ is strictly decreasing in $\left(x_1, \frac{n+1}{n-1}\right)$ and $g'(1) = 0$, we have that $g'(x_1) > 0$. But $g'(0) < 0$, so Bolzano's theorem yields that there exists $x_2 \in (0,x_1)$ such that $g'(x_2) = 0$. Given the monotonicity of $g'$ on that interval, $x_2$ must be the unique root there.

                    Summing up, we deduced that $g$ is decreasing on $[0,x_2]$ and on $\left[1, \frac{n+1}{n-1}\right]$; and increasing on $[x_2,1]$. Hence
                    \[
                        g(x_2) < g(1) = \dfrac{4}{n+1} \leq 1 = g(0).
                    \]
                    Again, Bolzano's intermediate values theorem yields that there exists a (unique) $x_0 \in (0, x_2)$ such that $g(x_0) = \frac{4}{n+1}$. Given the monotonicity of $g$, we will have that
                    \[
                        \max_{x_0 \leq x \leq \frac{n+1}{n-1}} \{ g(x) \} = \max \left \{ g(x_0), g(1) \right \} = \dfrac{4}{n+1}.
                    \]
            \end{proof}

            \begin{lem} \label{lem:AppendixLemma3} 
                Consider the polynomial
                    \[
                        g_\varepsilon(\zeta) = 14\varepsilon + 14\zeta - 14 \varepsilon \zeta^2 + \varepsilon^2 \zeta^3.
                    \]
                    Then there exists $0 < R_1 < 1 < R_2 < \left(\frac{14}{13}\right)^{13}$ and $\varepsilon_0 = \varepsilon_0 (R_1, R_2)$ such that
                    \begin{equation} \label{eqn:AppendixEquation1}
                         |g_\varepsilon(\zeta)|^2 \left(1 - \dfrac{13}{14} \sqrt[13]{|\zeta|} \right)^2 < 1 + 27 \left( \dfrac{13}{14} \right)^{26} \varepsilon^2, \quad R_1 \leq |\zeta| \leq R_2,
                    \end{equation}
                    when $0 < \varepsilon < \varepsilon_0$.
            \end{lem}

            \begin{proof}
                On the one hand, using the second order Taylor expansion at $x = 1$ of the function
                    \[
                        x \mapsto 14^2x^2 \left( 1 - \dfrac{13}{14} \sqrt[13]{x}\right)^2,
                    \]
                    we see that
                    \begin{equation} \label{eqn:aux5}
                        14^2 |\zeta|^2 \left(1 - \dfrac{13}{14} \sqrt[13]{|\zeta|}\right)^2 = 1 - \dfrac{14}{13} (1 - |\zeta|)^2 + O\big( \big|1 - |\zeta| \big|^3 \big).
                    \end{equation}
                    In particular, there exists $0 < R_1 < 1 < R_2 < \left( \frac{14}{13}\right)^{13}$ such that
                    \begin{equation} \label{eqn:aux7}
                        14^2 |\zeta|^2 \left(1 - \dfrac{13}{14} \sqrt[13]{|\zeta|}\right)^2 \leq 1 - (1 - |\zeta|)^2, \quad R_1 \leq |\zeta| \leq R_2.
                    \end{equation}
                    By shrinking $[R_1,R_2]$ if necessary, we may assume further that
                    \begin{eqnarray}
                        \left| \dfrac{1}{|\zeta|} - |\zeta|\right| & \leq & 3 \big| 1 - |\zeta| \big|, \label{eqn:alpha}\\
                        \left| \left( \dfrac{1}{|\zeta|} + |\zeta|\right)^2 + 4 - \dfrac{3|\zeta|^2}{7} \right| & \leq & 8, \label{eqn:beta}\\
                        \left| \left( \dfrac{1}{|\zeta|} + |\zeta| \right)^2 - \dfrac{|\zeta|^2}{7}\right| & \leq & \dfrac{27}{4} \left( \left( \dfrac{13}{14}\right)^{26} + \dfrac{3}{7} \right) \label{eqn:gamma}
                    \end{eqnarray}
                    and
                    \begin{equation}
                        \left| \dfrac{2|\zeta|^2}{7} - 4 + \dfrac{13}{14} \left( \dfrac{1}{|\zeta|} + 1\right)^2\right| \leq \dfrac{27}{4} \left( \left( \dfrac{13}{14}\right)^{26} - \dfrac{1}{7} \right) \label{eqn:delta}
                    \end{equation}
                    when $R_1 \leq |\zeta| \leq R_2$.
                
                On the other hand, if $\zeta=|\zeta|e^{it}$, then
                    \begin{eqnarray}
                    \dfrac{|g_\varepsilon(\zeta)|^2}{14^2 |\zeta|^2} & = & \left| 1+ \left(\dfrac{1}{\zeta} - \zeta \right) \varepsilon + \dfrac{\varepsilon^2 \zeta^2}{14}\right|^2 \nonumber \\
                            & = & 1 + \left(\left( \dfrac{1}{|\zeta|} + |\zeta| \right)^2 -\dfrac{|\zeta|^2 }{7} \right)\varepsilon^2 + 2 \left( \dfrac{1}{|\zeta|} - |\zeta|  \right) \varepsilon \cos t + \left(\dfrac{2|\zeta|^2}{7} - 4 \right)\varepsilon^2\cos^2 t  + O(\varepsilon^3), \label{eqn:aux6}
                    \end{eqnarray}
                    where the big-$O$ constant depends only on $R_1$ and $R_2$.

                We will split the proof in two parts:

                \begin{itemize}
                    \item If $\big|1 - |\zeta| \big| \geq 8 \varepsilon$, then
                        \[
                            \big( 1 - (1 - |\zeta|)^2\big) \left( 1 + 8\varepsilon^2 + 6 \varepsilon \big| 1 - |\zeta| \big|\right) \leq  (1 - 64 \varepsilon^2) \left( 1 + 56 \varepsilon^2 \right) < 1,
                        \]
                        because the polynomial
                        \[
                            x \mapsto (1 - x^2) \left( 1 + 8\varepsilon^2 + 6 \varepsilon x \right)
                        \]
                        has negative derivative in $8 \varepsilon < x$. Thus, we will have
                        \[
                            |g_\varepsilon(\zeta)|^2 \left(1 - \dfrac{13}{14} \sqrt[13]{|\zeta|} \right)^2 < 1 + O( \varepsilon^3),  \quad 8 \varepsilon \leq \big| 1 - |\zeta| \big|
                        \]
                        by \eqref{eqn:aux7}, \eqref{eqn:alpha}, \eqref{eqn:beta}  and \eqref{eqn:aux6}.

                    \item If $|1 - |\zeta|| < 8\varepsilon$, the identity 
                        \[
                            2 \left( \dfrac{1}{|\zeta|} - |\zeta|  \right) \varepsilon \cos t - \dfrac{14}{13} (1 - |\zeta|)^2 = \dfrac{13}{14}\left( \dfrac{1}{|\zeta|} + 1 \right)^2 \varepsilon^2 \cos^2 t - \left( \sqrt{\dfrac{14}{13}} (1 - |\zeta|) - \sqrt{\dfrac{13}{14}} \left( \dfrac{1}{|\zeta|} + 1 \right) \varepsilon \cos t\right)^2
                        \]
                        can be used to deduce via \eqref{eqn:gamma} and \eqref{eqn:delta} that
                        \begin{align*}
                            1 + \left(\left( \dfrac{1}{|\zeta|} + |\zeta| \right)^2 -\dfrac{|\zeta|^2 }{7} \right)\varepsilon^2 + &2 \left( \dfrac{1}{|\zeta|} - |\zeta|  \right) \varepsilon \cos t + \left(\dfrac{2|\zeta|^2}{7} - 4 \right)\varepsilon^2\cos^2 t - \dfrac{14}{13} (1 - |\zeta|)^2\\
                            & \leq 
                            1 + \left(\left( \dfrac{1}{|\zeta|} + |\zeta| \right)^2 -\dfrac{|\zeta|^2 }{7} \right)\varepsilon^2 + \left(\dfrac{2|\zeta|^2}{7} - 4  + \dfrac{13}{14}\left( \dfrac{1}{|\zeta|} + 1 \right)^2 \right)\varepsilon^2 \cos^2 t \\
                            &\leq  1 + \dfrac{27}{2} \left( \left( \dfrac{13}{14}\right)^{26} + \dfrac{1}{7} \right) \varepsilon^2.
                        \end{align*}

                        Combining this with \eqref{eqn:aux5} (recall that $\big|1 - |\zeta| \big| = O(\varepsilon)$ in this range) and \eqref{eqn:aux6}, we conclude that
                        \[
                            |g_\varepsilon(\zeta)|^2 \left(1 - \dfrac{13}{14} \sqrt[13]{|\zeta|} \right)^2 \leq 1 + \dfrac{27}{2} \left( \left( \dfrac{13}{14}\right)^{26} + \dfrac{1}{7} \right) \varepsilon^2 + O( \varepsilon^3). 
                        \]
                \end{itemize}

                In any case, since $\frac{27}{2} \left( \left( \frac{13}{14}\right)^{26} + \frac{1}{7} \right)<27\left( \frac{13}{14}\right)^{26} $, we may choose $\varepsilon_0$ depending only on $R_1$ and $R_2$ such that \eqref{eqn:AppendixEquation1} holds for all $0 < \varepsilon < \varepsilon_0$.
            \end{proof}
            
            \begin{rk} \label{rem:AppendixRemark2}
                Using the change of variables $\zeta = \left(\frac{14}{13} z^2 \right)^{13}$, the above lemma shows that the function $f_\varepsilon$ defined in \eqref{eqn:Countexample27} satisfies that
                \[
                    \max_{r_1 \leq |z| \leq r_2} \left \{ |f_\varepsilon' (z)|^2 (1 - |z|^2)^2 \right\} < 1 + 27 \left( \dfrac{13}{14} \right)^{26} \varepsilon^2
                \]
                if $0 < \varepsilon < \varepsilon_0$ for some $0 < r_1 < \sqrt{\frac{13}{14}} < r_2 < 1$.
            \end{rk}
            
            \begin{lem}\label{lem:auxiliary_max}
                Let $x\in \mathbb{R}$, $C\geq 0$ and $0<\beta<1 $. Then 
                \[
                \max\left\{0, x-C\right\}+\beta C\geq \beta \max \left\{ x,0\right\}.
                \]
            \end{lem}
            \begin{proof}
                If $x \leq 0$ the right hand side is $0$. If $0<x<C$, the inequality is reduced to $\beta C\geq \beta x$. Finally, if $x\geq C$, then $x-C+\beta C\geq \beta(x-C)+\beta C=\beta x $.
            \end{proof}
                
    \bibliographystyle{amsplain}
	\bibliography{bibliography}

\end{document}